\documentclass[11pt]{article}
\usepackage{geometry}                
\usepackage{graphicx}
\usepackage{amssymb}
\usepackage{epstopdf}
\def\'#1{\ifx#1i{\accent"13 \i}\else{\accent"13 #1}\fi}

\def\B{\mathcal{B}_3}

\def\qed{\hfill$\bullet$}
\def\eoc{\hfill$\circ$}

\newtheorem{theorem}{Theorem}
\newtheorem{lemma}{Lemma}
\newtheorem{corollary}{Corollary}[lemma]

\title{On panchromatic patterns}
\author{\sc Hortensia Galeana-S\'anchez \\ \tt hgaleana@math.unam.mx \and \sc Ricardo Strausz \\ \tt dino@math.unam.mx \and \\ --------------------------------------------------------- \\ \rm Instituto de Matem\'aticas \\ Universidad Nacional Aut\'onoma de M\'exico \\ Ciudad Universitaria, 04510, M\'exico D.F. \\ ---------------------------------------------------------}

\begin{document}

\maketitle

\begin{abstract}
Given $D$ and $H$ two digraphs, $D$ is $H$-coloured iff the arcs of $D$ are coloured with the vertices of $H$. After defining what do we mean by an $H$-walk in the coloured $D$, we characterise those $H$, which we call {\sl panchromatic patterns\/}, for which all $D$ and all $H$-colourings of $D$ admit a kernel by $H$-walks. This solves a problem of Arpin and Linek from 2007 \cite{AL}.
\end{abstract}

\section{Introduction}

Tha notion of {\sl kernel\/}, as introduced by Von Neumann in 1944 \cite{VN}, has been generalised in several directions, see for example \cite{H,LS,R}. Here, we study yet another generalisation introduced by Sands et al. \cite{SSW} which arises naturally when the arcs of a digraph are coloured (see also \cite{GS}); such a generalisation was deeply studied too by Arpin and Linek \cite{AL}.

Following Arpin and Linek \cite{AL}, let $\B$ be the class of digraphs $H$ with the property that for every digraph $D$ and every $H$-colouring of $D$ (to be defined) there exists an $H$-kernel by walks in $D$. We call those digraphs in $\B$ {\sl panchromatic patterns\/}.

The aim of this paper is to characterise panchromatic patterns. For, we first prove a technical lemma (see Section~1) that allow us to add, under controlled circumstances, an arc to a panchromatic pattern preserving this property. This allow us to settle a question raised by Arpin and Linek \cite{AL} which characterises all panchromatic patterns of order 3.

Then, we introduce the notion of a {\it bicomplete\/} digraph and prove (see Lemma~\ref{bicompletas}) the sufficiency for such digraphs to be panchromatic patterns. The rest of the paper, after some preliminaries in Section~2, is devoted to prove the necessity of such a property to be a panchromatic pattern and settle the desired characterisation.

\section{Preliminaries}

By a {\it digraph\/} $D=(V,A)$ we mean a finite non-empty set of vertices $V$ and a set of (directed) arcs $A\subseteq V\times V$. Given another digraph $H$, by an {\it $H$ colouring of $D$\/}, we mean a map $\varsigma\colon A(D)\to V(H)$ from the arcs of $D$ to the vertices of $H$ --- we think on the vertices of $H$ as colours assigned to the arcs of $D$, hence the name. Given such a colouring, a walk $W=x_0,x_1,\dots,x_k$ in $D$ is called an {\it $H$-walk\/} if $\varsigma(W)=\varsigma(x_0,x_1),\varsigma(x_1,x_2),\dots,\varsigma(x_{k-1},x_k)$ is a walk in $H$. A subset $K\subset V(D)$ is called an {\it $H$-kernel\/} if it is both, $H$-independent and $H$-absorbent; viz., there are no $H$-walks between any pair of different vertices in it, and given any vertex out of it, there is an $H$-walk into such a subset. 

Suppose that $H$ is a looped digraph with the following properties:
\begin{enumerate}
\item $V(H) = \bigsqcup_{i=1}^n C_i$,
\item $(x,y)\in A(H)$, whenever $x\not=y$ and $x,y\in C_i$ for some $i$, and
\item $C_i\times C_j\subseteq A(H)$ whenever $i\not=j$ and $C_i\times C_j\cap A(H)\not=\emptyset$ or $i=j$ and $(x,x)\in A(H)$ for some $x\in C_i$.
\end{enumerate}
Then, the digraph $H'$ with vertices $V(H')=\{C_1,\dots,C_n\}$ and arcs $(C_i,C_j)\in A(H')$ if and only if $(C_i\times C_j)\cap A(H)\not=\emptyset$ is called a {\it contraction\/} of $H$ and $H$ is an {\it expansion\/} of $H'$.

A digraph $H$ is called a {\it panchromatic pattern\/} if given any digraph $D$ and any $H$-colouring of $D$, we can find an $H$-kernel of $D$. The aim of this paper is to characterise the class $\B$ of panchromatic patterns. For, we will use the following results of Arpin and Linek \cite{AL}, without further reference:

\begin{lemma}
\label{inducidas}\ \\
\begin{enumerate}
\item If $H\in\B$, then every vertex of $H$ has a loop (is looped),
\item If $H\in\B$, and $H'$ is an induced subdigraph of $H$, then $H'\in\B$,
\item Let $H'$ be a contraction of $H$. $H\in\B$ if and only if $H'\in\B$.
\end{enumerate}

\end{lemma}

\begin{lemma}
\label{caminos}
Let $W=x_0,x_1,\dots,x_k$ be a walk in $H$ such that
\begin{enumerate}
\item for all $x_j\in W$, with $0\leq j\leq k-1$, there is a colour $c_j\in V(H)$ such that $(x_j,c_j)\not\in A(H)$,
\item $(x_k,x_0)\not\in A(H)$.
\end{enumerate}
Then, $H\not\in\B$.
\end{lemma}

\begin{figure}[htbp] 
   \centering
   \includegraphics[width=4in]{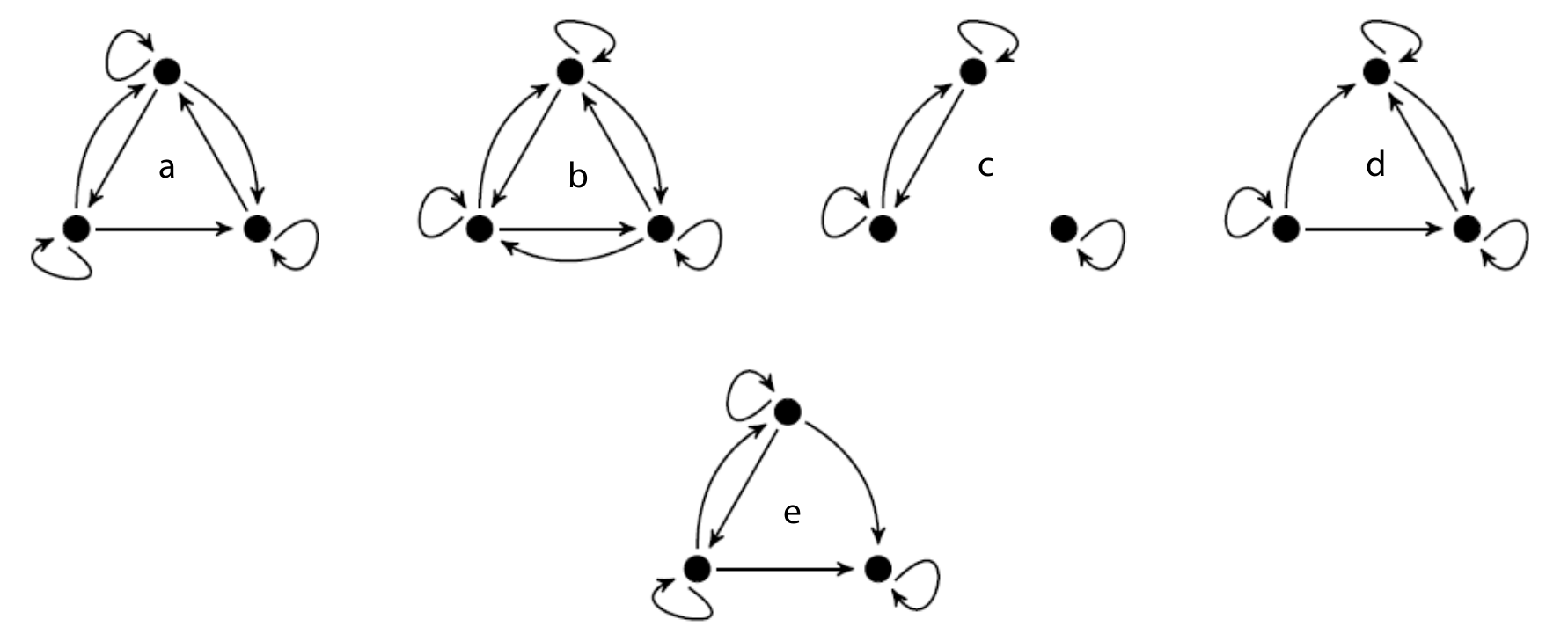} 
   \caption{Panchromatic patterns of order 3}
   \label{enb3}
\end{figure}

\begin{lemma}
The digraphs depicted in Figure~\ref{enb3} are all elements of $\B$.
\end{lemma}

\begin{figure}[htbp] 
   \centering
   \includegraphics[width=4in]{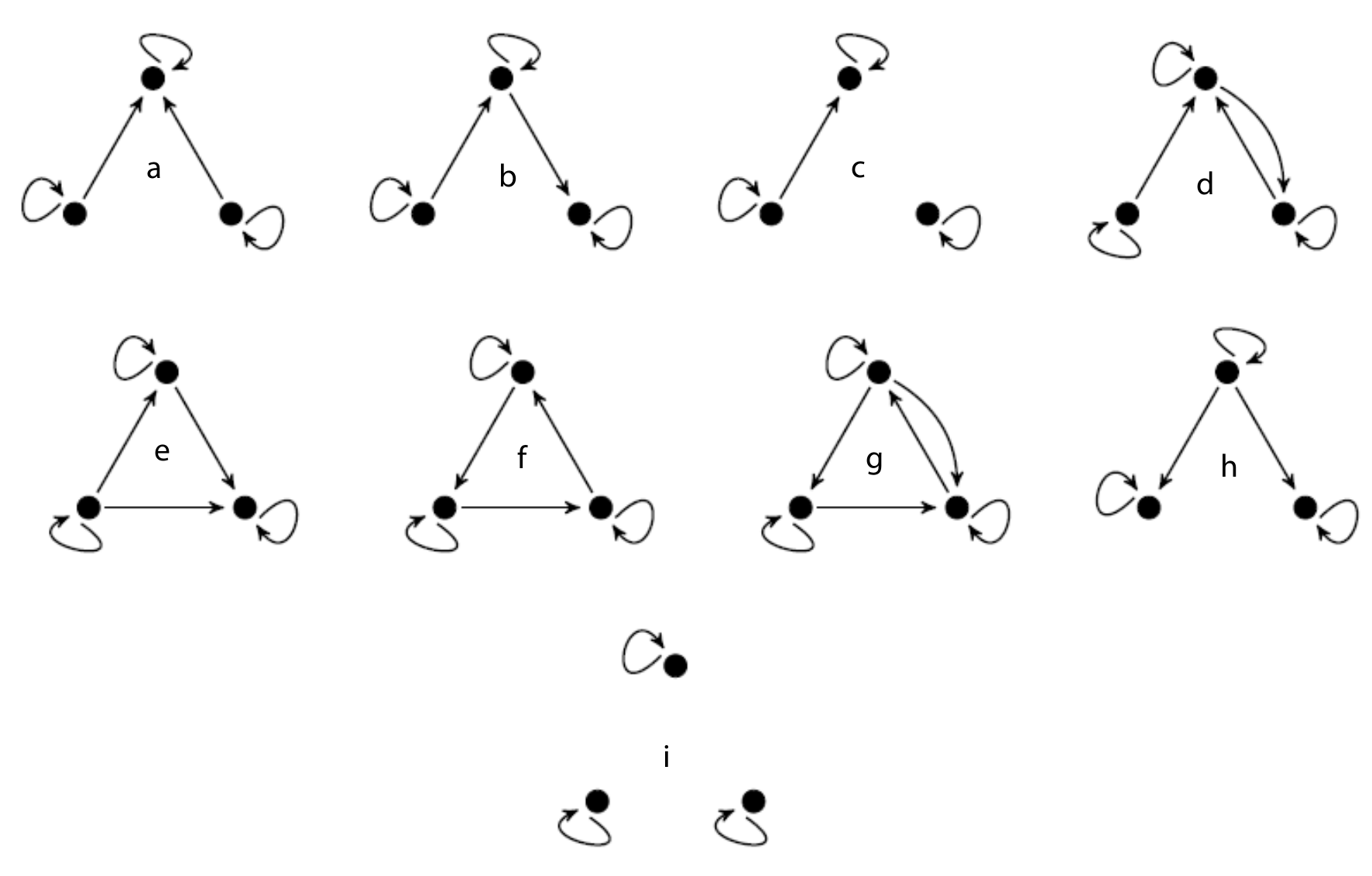} 
   \caption{Some non-panchromatic patterns of order 3}
   \label{noenb3}
\end{figure}

\begin{lemma}
Non of the digraphs depicted in Figure~\ref{noenb3} are elements of $\B$.
\end{lemma}

\begin{lemma}
\label{impares}
If $H\in\B$, then $H$ does not contain odd directed cycles in its complement.
\end{lemma}

\section{Main Lemmas}

\begin{lemma}
\label{p2}
Let $H$ be a digraph, with all its vertices looped, and $a=(u,v)\not\in A(H)$ a pair of vertices at distance 2 from $u$ to $v$. If $H\cup a\not\in\B$, then $H\not\in\B$.
\end{lemma}

\noindent
{\bf Proof.} Let $H$ be a looped digraph and $u,z,v$ a path in $H$ from $u$ to $v$, two non adjacent vertices of $H$. Let $H' = H\cup(u,v)$. 

\begin{figure}[htbp] 
   \centering
   \includegraphics[width=3in]{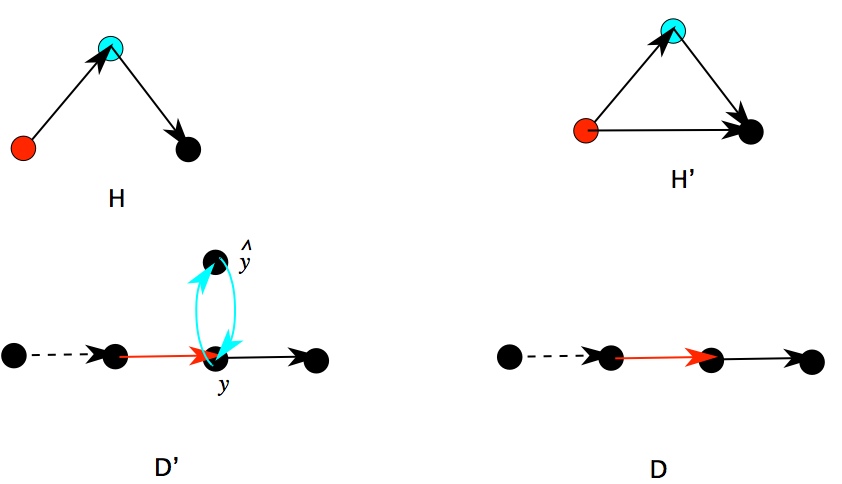} 
   \caption{The construction}
   \label{construccion}
\end{figure}

First we will show that for every $D$ and every $H'$-colouring of $D$, there exists a digraph $D'$ and an $H$-colouring of $D'$, such that $D$ admits an $H'$-kernel by walks if (and only if) $D'$ admits an $H$-kernel by walks. For, let $D'$ be constructed from $D$ as follows: for each walk $(x,y,z)$ in $D$, with arcs coloured by $u$ and $v$, in that order, we add a new vertex $\hat y$ in $D'$ and the symmetric arrows $(y,\hat y)$ and $(\hat y,y)$, each of them coloured with $z$. We denote by $\hat Y = \{\hat y\}$ the set of all those new vertices in $D'$ and by $Y$ their corresponding neighbours in $D$ (see Figure~\ref{construccion}).

Let $K'\subset V(D')$ be an $H$-kernel by walks. We can construct an $H'$-kernel by walks in $D$ with the following set
	$$K:=K'\cup\{y\in V(D) : \hat y\in K'\}\setminus\hat Y.$$

\begin{figure}[htbp] 
   \centering
   \includegraphics[width=1.3in]{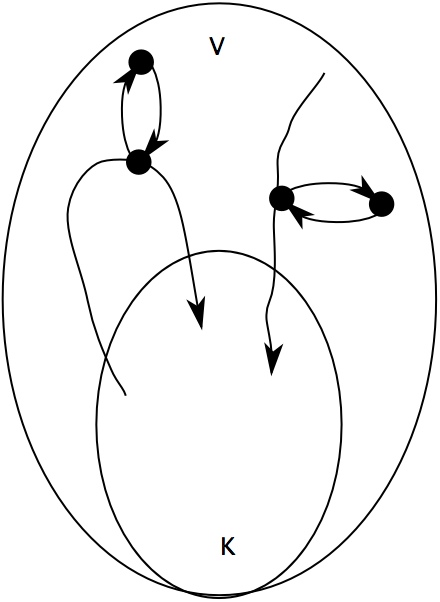} 
   \caption{The proof}
   \label{prueba}
\end{figure}

To see the independence of $K$, consider a couple of vertices $\alpha,\beta\in K$. If there is an $H'$-walk from $\alpha$ to $\beta$ and such a walk never uses a vertex in $Y$, such a walk is an $H$-walk and it would contradict the independence of $K'$ in $D'$. So, let us suppose that such walk passes through a vertex $y\in Y$. Then, we can construct an $H$-walk in $D'$ adding the neighbour $\hat y$ of $y$ and the arrows $(y,\hat y)$ and $(\hat y,y)$ in the position of the vertex $y$ (see Figure~\ref{prueba}).

We have to be careful and note that, if $\alpha$ or $\beta$ are vertices of $Y$, we need only to add the arcs $(\hat y,y)$ or $(y,\hat y)$, respectively, at the beginning or the end of the $\alpha\beta$ walk to obtain an $H$-walk in $D'$ and contradicting the independence of $K'$.

Analogously, we can show the absorbency of $K$ by omitting the occurrences of vertices in $\hat Y$ in the $H$-walks of $D'$ and using the $H$-absorbency of $K'$. 

Finally, seeking for a contradiction, suppose that $H\in\B$. Since $H'\not\in\B$, then there exists a $D$ and an $H'$-colouring of $D$ which do not admit an $H'$-kernel by walks. Let $D'$ be constructed and coloured as before. Since $H\in\B$ then $D'$ admits an $H$-kernel by walks. Therefore, by the previous argument, $D$ would admit an $H'$-kernel by walks contradicting the hypothesis of the theorem, and concluding the proof. \qed

\begin{figure}[htbp] 
   \centering
   \includegraphics[width=2in]{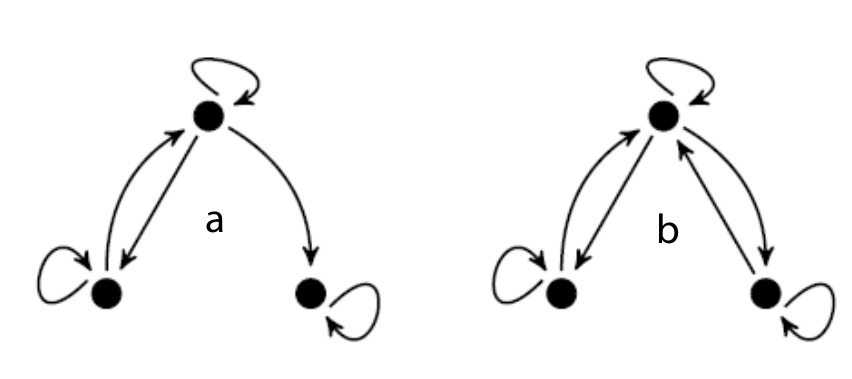} 
   \caption{Patterns previously not known if panchromatic}
   \label{nuevas}
\end{figure}

As a consequence of this lemma, we can now decide that the following two digraphs in Figure~\ref{nuevas} are not in $\B$, which was previously unknown (see again \cite{AL}).

\begin{corollary}
The patterns in Figure~\ref{nuevas} are not panchromatic.
\end{corollary}

\begin{figure}[htbp] 
   \centering
   \includegraphics[width=1in]{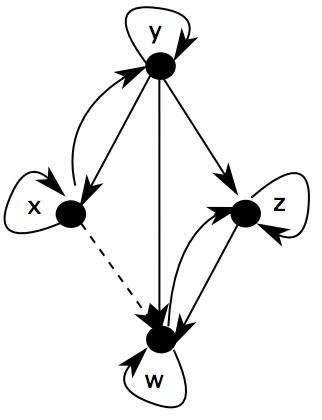} 
   \caption{Extension of Figure \ref{nuevas}.a}
   \label{auxiliara}
\end{figure}

\noindent
{\bf Proof.}  We first concentrate on Figure~\ref{nuevas}.a. For, due to Lemma~\ref{inducidas}.3, it is enough to show that the solid arrows in Figure~\ref{auxiliara} induce a non-panchromatic digraph. Consider such a digraph and add the doted arc $(x,w)$. Due to lemma~\ref{p2}, it is enough now to show that this new digraph is not a panchromatic pattern; this last follows from the fact that the subdigraph induced by $x,w$ and $z$ is not a panchromatic pattern (see Figure~\ref{noenb3}.d).

\begin{figure}[htbp] 
   \centering
   \includegraphics[width=1in]{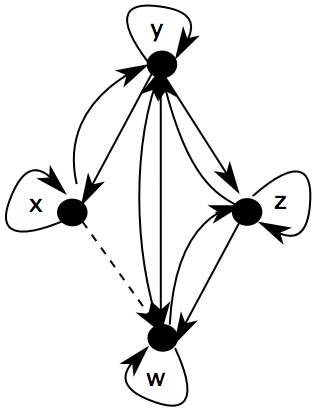} 
   \caption{Extension of Figure \ref{nuevas}.b}
   \label{auxiliarb}
\end{figure}

Analogously, to show that Figure~\ref{nuevas}.b is not a panchromatic pattern, we extend it with a vertex $w$ (see figure~\ref{auxiliarb}), add the dotted arc, and find the subdigraph induced by $x,w$ and $z$, which we know is not in $\B$.

\qed

\bigskip

We say that a looped digraph $H$ is {\sl bicomplete} if it has the following three properties:
\begin{enumerate}
\item $V(H) = X\sqcup Y$,
\item $H[X]$ and $H[Y]$ are complete digraphs, and
\item $\forall x\in X$ and $y\in Y$ : $(y,x)\in A(H)$.
\end{enumerate}
Note that in a bicomplete digraph some of the $(x,y)$ arcs may belong to $A(H)$.

\begin{lemma}
\label{bicompletas}
If $H$ is bicomplete, then $H\in\B$.
\end{lemma}

\noindent
{\bf Proof.} Let $H$ be a bicomplete digraph. On the one hand, let us suppose that there is no arc $(x,y)\in A(H)$, with $x\in X$ and $y\in Y$. Consider the digraph 
	$$H':=\left(V(H)\cup\{z\},A(H)\cup\{(x,z),(z,x),(y,z),(z,y),(z,z) : x\in X, y\in Y\}\right).$$
It is easy to se that $H'$ can be contracted to the digraph in Figure \ref{enb3}.a, which is in $\B$, and therefore by the Lemma~\ref{inducidas}.3, $H'$ and all its induced subdigraphs (Lemma~\ref{inducidas}.2) belongs to $\B$ --- in particular, $H\in\B$.

On the other hand, consider a new pair $a=(x,y)$. We can add the arc $a$ to $H$ preserving its panchromaticity (i.e., $H\cup a\in\B$). For, observe that there is a path of length 2 $(x,z,y)$ in $H'$ and, by the previous Lemma~\ref{p2}, we guaranty that $H'\cup(x,y)\in\B$. Since $H\cup a$ is an induced subdigraph of $H'\cup a$, by Lemma~\ref{inducidas}.2 we guarantee its panchromaticity.

Thus, we can recursively add arcs from $X$ to $Y$ and in each step we preserve the panchromaticity of the pattern.\qed

\section{Main Theorem}

Due to Sands et al. \cite{SSW} we know that the digraph $2K_1$ consisting in two looped vertices belongs to $\B$; moreover, due to Lemma~\ref{inducidas}.3, any expansion of $2K_1$ belongs to $\B$. Furthermore, due to lemma~\ref{bicompletas} we know that every bicomplete digraph is also in $\B$. We will now show that these are all digraphs in $\B$.

\begin{theorem}
$H$ is a panchromatic pattern if and only if $H$ is bicomplete or $H$ can be contracted to $2K_1$.
\end{theorem}

\noindent
{\bf Proof.}  The sufficiency is the content of Lemma~\ref{bicompletas} and Sauer et al. \cite{SSW}, respectively. For the necessity, let us suppose that $H\in\B$ is a panchromatic pattern and that it is not an expansion of $2K_1$. Let us denote by $G=H^c$ the complementary digraph of $H$.

{\bf Claim.} {\sl If $(u,v)\in A(G)$ is an asymmetric arc, then $d_G^+(v)=0$.\/}

For, suppose that there exist a $y\not=u$ such that $a=(v,y)\in A(G)$. If $a$ is an asymmetric arc (i.e., $(y,v)\not\in A(G)$), then, depending on the relationship between $u$ and $y$, one of the digraphs in Figure~\ref{noenb3}.b,e,f,g appears as an induced subdigraph of $H$, contradicting its panchromaticity. If $a$ is symmetric (i.e., $(y,v)\in A(G)$), depending on the relationship between $u$ and $y$, one of the digraphs in Figure~\ref{noenb3}.a,b,c,d appears as an induced subdigraph of H, contradicting again its panchromaticity. In any case, $a$ cannot belong to $G$, concluding the proof of the claim. \eoc

As an immediate consequence of this claim we conclude that: {\sl if $H\in\B$ then every cycle of length at least 4 is symmetric}. 

Recall that Arpin and Linek \cite{AL} proved that no digraph in $\B$ has odd {\sl directed\/} cycles in its complement.

{\bf Claim.} {\sl The underling graph of $G$ has no odd cycles; i.e., $G$ is bipartite.}

Searching for a contradiction, let us suppose that the underling graph of $G$ has an odd cycle. Such a cycle cannot be symmetric since it would contain a directed odd cycle, contradicting Lemma~\ref{impares}, therefore it has an asymmetric arc. Due to the previous claim, it is easy to see that Figure~\ref{nuevas}.a, without its loops, most be part of the supposed cycle. From here we conclude that one of the digraphs in Figure~\ref{nuevas}.a or Figure~\ref{noenb3}.b,c,h most be induced subdigraphs of $H$ contradicting its panchromaticity. \eoc

From here, we have two cases to analyse; namely, if $G$ has directed cycles of length at least 4, or not.

\bigskip
\noindent
{\bf Case 1.} {\sl If $G$ contains a directed cycle of length at least 4, then $G$ is a bipartite complete digraph, and therefore $H$ is an extension of $2K_1$. \/}

For, recall that every directed cycle of length at least 4 is symmetric. Therefore every non-trivial strongly connected component of $G$ is symmetric (since in a strongly connected digraph each arc is in a directed cycle).

{\bf Claim.} {\sl If $S$ is a strongly connected component of $G$ with a directed cycle of length at least 4, then $S$ is a bipartite complete digraph (i.e., all arcs are symmetric and it has all arcs between two independent sets of vertices).\/}

By the previous claim, $G$ is bipartite and therefore the induced subgraph of $S$ is bipartite too. Let $V(S)=A\sqcup B$ be the bipartition of $S$. We first show that $S$ has a cycle of length exactly 4. Let $\gamma$ be a directed girth of $S$. Searching for a contradiction, suppose $|\gamma|>4$; since there are no odd cycles, we have that $|\gamma|\geq6$. This induces a $P_4=(x,w,v,u)$ subgraph in $G$. Then, we have the path $(w,u,x,v)$ induced in $H$ and by Lemma~\ref{caminos} we have that $H\not\in\B$. Thus, let $(u,v,w,x)$ be the cycle of length 4 guaranteed by the previous argument. We further suppose that $u$ and $w$ are in $A$.

Observe that every vertex in $B$ is adjacent to $u$ (by a symmetric arc); for, if there exists a vertex $z$ in $B$ non adjacent to $u$, then the induced subgraph by $u$, $z$ and $v$ is isomorphic to Figure~\ref{nuevas}.b, contradicting the panchromaticity of $H$. Analogously $w$ is adjacent to each vertex in $B$; moreover, $v$ and $x$ are adjacent to all vertices in $A$. Furthermore, each vertex $y\in A\setminus\{u,w\}$ is adjacent to all vertices in $B$; for, observe that if there is a vertex $z\in B\setminus\{v,x\}$ non adjacent to $y$, then the subdigraph induced by $w$, $y$ and $z$ is isomorphic to Figure~\ref{nuevas}.b, which contradicts the panchromaticity of $H$, concluding the proof of the claim.\eoc

{\bf Claim.} {\sl Every connected component of $G$ is strongly connected.\/}

For, let $S$ be as before with its 4-cycle $(u,v,w,x)$ and bipartition $V(S)=A\sqcup B$, and suppose that the connected component of $S$, in the week sense, is bigger. Then either an asymmetric arc goes into $S$ from its complement, or there is an asymmetric arc from $S$ to its complement. Due to our first claim, there cannot be an ``incoming'' arc $(y,s)$, with $s\in S$ and $y\in G\setminus S$, so let us suppose there is an ``outgoing'' arc $(s,y)$. With out loose of generality, suppose that $s\in B$. Then, the subdigraph induced by $w$, $s$ and $y$ in $H$, depending in the relationship between $w$ and $y$,  is isomorphic to one of the Figures~\ref{noenb3}.b,c,h or Figure~\ref{nuevas}.a, contradicting the panchromaticity of $H$ and concluding the claim.\eoc

{\bf Claim.} {\sl The underlying graph of $G$ is connected.\/}

Let $S$ be as before, suppose there is another component of $G$, and let $y$ be a vertex in $G\setminus S$. Then the subdigraph of $H$ induced by $w$, $x$ and $y$ is isomorphic to Figure~\ref{nuevas}.b, which contradicts the panchromaticity of $H$.\eoc

Therefore, we have that $G$ is isomorphic to $S$ which we have shown to be a complete bipartite and $H$ is an extension of $2K_1$ concluding case~1.

\bigskip
\noindent
{\bf Case 2.} {\sl If $G$ does not have a directed cycle of length at least 4, then $H$ is a bicomplete digraph.\/}

Recall that we had proved that the underling graph of $G$ is bipartite; we will work with its bipartition $V(G)=A\sqcup B$.

{\bf Claim.} {\sl If $G$ contain cycles (viz., symmetric arcs), then all of them pass through a single vertex $x\in V(G)$.\/}

If $G$ is acyclic, we are done; so, let us suppose that $G$ has a symmetric arc $\{u,v\}$, where $u\in A$ and $v\in B$. We will show that either every cycle of $G$ contain $u$ or all of them contain $v$. First of all, suppose there is a cycle (of length 2 --- or a symmetric arc, if you will) that does not contain either $u$, nor $v$; call such an arc $\{z,w\}$ where $z\in A$ and $w\in B$. Then the symmetric arc $\{v,z\}$ must exist since the complement of $v,w$ and $z$ must be in Figure~\ref{enb3}. Analogously, the pair $\{u,w\}$ form a symmetric arc; therefore we have the cycle $(u,v,z,w)$ which contradict the fact that every cycle is of length 2.

Thus, we can suppose now that every cycle either contains $u$, or it contains $v$. Suppose there are cycles $\{z,v\}$ and $\{u,w\}$ in $G$. Since we don't have cycles of length 4, $z$ and $w$ must be independent, therefore we have that the complement of $u,z$ and $w$ is Figure~\ref{nuevas}.b, contradicting the panchromaticity of $H$ and concluding the propf of the claim. \eoc

{\bf Claim.} {\sl There is a partition of $V(G\setminus x)=A\sqcup B$ (possibly degenerated; i.e., $B=\emptyset$), where $x$ is that vertex contained in all symmetric arcs, such that all arcs go from $A$ to $B$.\/}

For, simply observe that $G\setminus x$ is acyclic and that, as proved earlier, the final vertex of every asymmetric arc has out-degree zero; so, let $A$ be the set of initial vertices of all arcs, and $B$ its complement. \eoc

We end the proof showing that such a vertex $x\in V(G)$, indeed, does not exist.

{\bf Claim.} {\sl If $H$ is not an extension of $2K_1$, then the digraph $G$ is acyclic.\/}

For, suppose we have such a vertex $x\in V(G)$ and the partition $V(G\setminus x)=A\sqcup B$, where all arcs go from $A$ to $B$. By definition of $x$, there is a symmetric arc $\{x,y\}$ in $G$. If $y$ is a vertex in $A$, then there are symmetric arcs $\{x,a\}$ for all vertices in $A$ since otherwise the complement of $x,y$ and the nonadjacent $a\in A$ would induce a digraph isomorphic to Figure~\ref{noenb3}.d or one of the two digraphs in Figure~\ref{nuevas}. Furtheremore, $B$ has to be empty since otherwise $x,a$ and $b$ would induce a subgraph isomorphic to Figure~\ref{noenb3}.b,c,h or Figure~\ref{nuevas}.a,b. So, $H$ would be an isolated vertex union a complete digraph, which is an extension of $2K_1$. Analogously, there are no symmetric arcs from $x$ to $B$.\eoc

Finally, since $G$ is bipartite and acyclic, and all out degrees of final vertices are 0, then either $H$ is bicomplete or an extension of $2K_1$, which completes the proof.\qed


\end{document}